\documentclass[reqno]{amsart}

\pdfoutput = 1

\usepackage[english]{babel}
\usepackage{amsmath,amssymb,amsthm}
\usepackage{amsfonts}
\usepackage{array,dcolumn}
\usepackage{graphicx}

\usepackage[T1]{fontenc}
\usepackage[utf8]{inputenc}
\usepackage[activate={true,nocompatibility},spacing,kerning]{microtype}
\usepackage[hyperpageref]{backref}
\usepackage[pdftex,
            colorlinks,
            hyperfootnotes=false,
            pdffitwindow=true,
            plainpages=false,
            pdfpagelabels=true,
            pdfpagemode=UseOutlines,
            pdfpagelayout=SinglePage,
            pdftitle={Convergence Rates of Gauss and Clenshaw--Curtis Quadrature},
            pdfauthor={Folkmar Bornemann, Technische Universität München},
            pagebackref,
            hyperindex,
            filecolor=purple,
            urlcolor=purple,
            citecolor=blue,
            linkcolor=blue]{hyperref}
\usepackage{myharvard}
\usepackage{bm}
\usepackage[sc,osf]{mathpazo}
\microtypecontext{spacing=nonfrench}

\graphicspath{{plots/}}

\renewcommand{\leq}{\leqslant}
\renewcommand{\geq}{\geqslant}

\newcommand{\R}{\mathbb{R}}

\DeclareMathOperator*\primesum{{\sum\nolimits\!'}}

\newtheorem*{thm}{Theorem}

\begin{document}

\renewcommand{\figurename}{\footnotesize {\sc Figure\@}}

\title[Convergence Rates of Gauss and Clenshaw--Curtis Quadrature]{On the Convergence Rates of Gauss and Clenshaw--Curtis Quadrature for Functions of Limited Regularity}

\author{Shuhuang Xiang}
\address{Dept. of Appl. Math., Central South University, Changsha, Hunan 410083, P. R. China}
\email{xiangsh@mail.csu.edu.cn}
\author{Folkmar Bornemann}
\address{Zentrum Mathematik -- M3, Technische Universität München,
         80290~München, Germany}
\email{bornemann@tum.de}

\begin{abstract}
We study the optimal general rate of convergence of the $n$-point quadrature rules of Gauss and Clenshaw--Curtis when applied
to functions of limited regularity: if the Chebyshev coefficients decay at a rate $O(n^{-s-1})$ for some $s>0$,
Clenshaw--Curtis and Gauss quadrature inherit exactly this rate. The proof (for Gauss, if $0<s<2$, there is numerical evidence only) is based on work of Curtis, Johnson, Riess, and Rabinowitz from the early 1970s and on a refined estimate for Gauss
quadrature applied to Chebyshev polynomials due to \citeasnoun{MR1345417}. 
The convergence rate  of both quadrature rules is up to one power of $n$ better than polynomial best approximation; hence, the classical
proof strategy that bounds the error of a quadrature rule with positive weights 
by polynomial best approximation is doomed to fail in establishing the optimal rate.
\end{abstract}

\maketitle

\section{Introduction}\label{sect:intro}

Though Clenshaw--Curtis and Gauss quadrature are classical topics in numerical analysis, it is quite hard to track down a theorem 
that would establish the \emph{optimal} rate of the error $E_n(f)$ of the
$n$-point rules for functions $f : [-1,1] \to \R$ of limited regularity. 
Here,
regularity is most conveniently measured\footnote{Some ways to determine $s$ are discussed in §\ref{sect:Xs}.} by the exponent $s > 0$ of a decay rate $a_m = O(m^{-s-1})$ 
of the coefficients $a_m$ of the expansion
\[
f(x) = \primesum_{m=0~\,}^{\infty~\,} \!\! a_m T_m(x)
\]
in terms of the Chebyshev polynomials $T_m(x)$ of the first kind of degree $m$; the prime indicates that the first term
is to be halved. We say that such a function $f$ is of class~$X^{s}$ and
claim that the error of both quadrature rules inherits exactly this rate:
\begin{equation}\label{eq:claim}
E_n(f) = O(n^{-s-1}).
\end{equation}
As noted by \citeasnoun[p.~893]{MR2600548}, the case $s=1$ can be found explicitly in the classical literature (we denote by $E_n^C(f)$ the quadrature error
of Clenshaw--Curtis and by $E_n^G(f)$ that of Gauss): if $f \in X^1$,
\begin{itemize}
\item \citeasnoun{MR0305555} proved  $E_n^C(f) = O(n^{-2})$;\\*[-3.5mm]
\item \citeasnoun[§4.8]{MR760629} gave a sketch that $E^G_n(f) = O(n^{-2})$.
\end{itemize}   
It is a fairly straightforward exercise, however, to extend the approach taken by these authors to the case of general $s>0$: an approach that starts from the bound
\begin{equation}\label{eq:series}
|E_n(f)| \leq \sum_{m=n}^{\infty} |a_m|\cdot |E_n(T_m)|.
\end{equation}
By using aliasing of under-sampled trigonometric polynomials, \citeasnoun{MR0305555} and \citeasnoun{MR0298934} showed, for Clenshaw--Curtis and Gauss
quadrature, that $E_n(T_m)$ is, up to some remainder, periodic in $m$ with a period of $O(n)$ and 
an \emph{average} modulus of $O(n^{-1})$. Hence, provided the remainder can effectively be controlled, one 
would read off 
the rate (\ref{eq:claim}). If it were not for this proviso, the story could end here; but the precise state 
of affairs differs considerably:
\begin{itemize}
\item For Clenshaw--Curtis quadrature, the remainder is a term of higher order, indeed; its effective control established by
\citeasnoun{MR0305555} for $s=1$ easily carries over to $s>0$; see §\ref{sect:cc} of this paper.\\*[-3.5mm]
\item For Gauss quadrature, the sketch given by \citeasnoun[§4.8]{MR760629} \emph{neglects} the remainder. 
Since it is not of strictly higher order, the remainder is much harder to control: aliasing holds asymptotically up to $m=o(n^{3/2})$ only; for larger $m$, phase errors of order~$O(1)$ enter. 
\end{itemize} 
Accordingly, to rigorously deal with Gauss quadrature, we split 
(\ref{eq:series}) after the first $O(n^{3/2})$ terms; the tail is then easily estimated by the decay of the coefficients and
a simple uniform bound of $E_n(T_m)$; see §\ref{sect:GI}. Using the estimate of the remainder given by \citeasnoun{MR0298934}, we are
able to prove the rate (\ref{eq:claim}) up to a factor $\log n$ for $s\geq 2$, whereas the case $0<s<2$ yields a 
suboptimal $O(n^{-3s/2})$ bound. 
Using a refinement of the Curtis--Rabinowitz estimate due to \citeasnoun{MR1345417}, \citeasnoun{Xiang} has recently eliminated the logarithmic factor for $s\geq 2$ (there is, still, no improvement
in the case $0<s<2$); see §\ref{sect:GII}. 

Summarizing, we have proved (\ref{eq:claim}) for all cases except for Gauss with $0< s < 2$:

\begin{thm} If $f \in X^s$, the error of $n$-point Clenshaw--Curtis quadrature and, for $s\geq 2$, also that of Gauss quadrature have the rate $O(n^{-s-1})$.
For $0<s<2$, the Gauss quadrature error is (at most) of size $O(n^{-3s/2})$.
\end{thm}

\begin{figure}[tbp]
\begin{center}
\begin{minipage}{0.495\textwidth}
{\includegraphics[width=\textwidth]{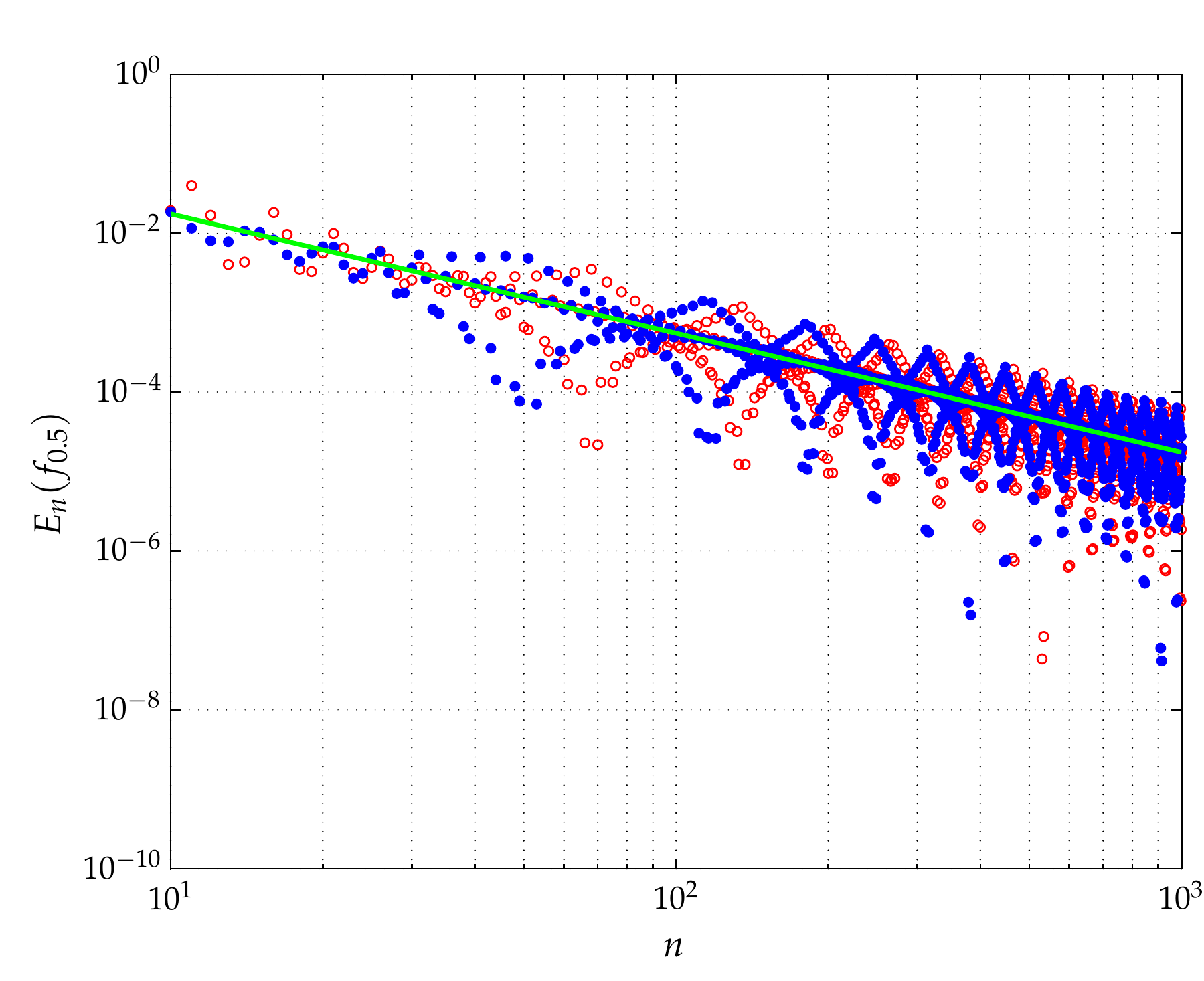}}
\end{minipage}
\hfil
\begin{minipage}{0.495\textwidth}
{\includegraphics[width=0.985\textwidth]{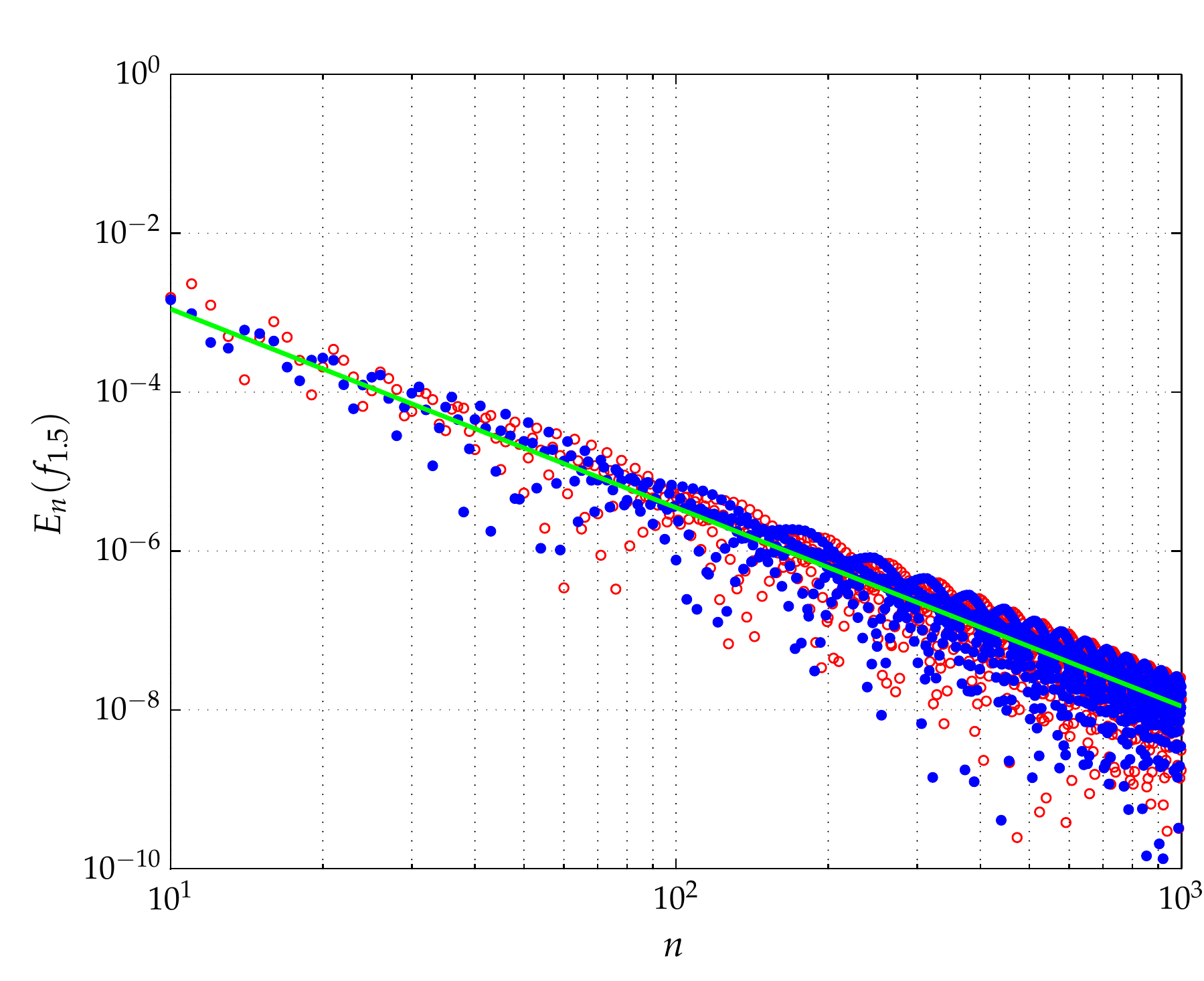}}
\end{minipage}
\end{center}~\\*[-5mm]
\caption{\footnotesize Numerical evidence that $n$-point Gauss quadrature has an $O(n^{-s-1})$ error rate for integrating the functions $f_s(x) = |x-0.3|^s$ (left: $s = 0.5$, right: $s=1.5$) on the interval $(-1,1)$: $E_n^G(f_s)$ (dots), $E_n^C(f_s)$ (circles), $c_s n^{-s-1}$ (solid line).}
\label{fig:1}
\end{figure}

Numerical experiments with $f_s(x)=|x-0.3|^s$, which is of class $X^s$ (see §\ref{sect:Xs}), and various $0<s<2$ (as in Fig.~\ref{fig:1}) has led us to the conjecture that 
Gauss quadrature enjoys the same $O(n^{-s-1})$ error rate as Clenshaw--Curtis also for $0<s<2$ in general. We remark that these experiments also show that the $O(n^{-s-1})$ error rate cannot be improved for any of the two quadrature rules.

\subsubsection*{Quadrature vs. best approximation}

In his detailed study of the almost equal numerical performance of the quadrature rules of Gauss and Clenshaw--Curtis
for functions of various regularity types, \citeasnoun{MR2403058} 
proved a suboptimal~$O(n^{-s})$ bound for functions $f\in X^s$. In the
Gauss case 
 he based his rate estimate on the classical bound $|E_n^G(f)|\leq 4E_{2n+1}^*(f)$ \citeaffixed[p.~333]{MR760629}{see, e.g.,} where $E_n^*(f)$ denotes the error of best approximation
by polynomials of degree $n$; if $f \in X^s$ this allows the straightforward estimate \citeaffixed[Thm.~3.3]{MR1060735}{see, e.g.,}
\[
E^*_n(f) \leq \sum_{m=n+1}^{\infty} |a_m| = O(n^{-s}).
\]
In the case $f(x)=|x|$ (which is of class $X^1$) the estimate is \emph{sharp}, since it is known by a theorem of Bernstein that \cite[Eq.~(1.18)]{MR891763}
\[
\lim_{n\to\infty} nE_n^*(|x|) = 0.2801694990\ldots\,.
\] 
Hence,  Clenshaw--Curtis and Gauss quadrature converge with a rate that is typically one power of 
$n$ better than the one of polynomial best approximation. 

\section{Functions of class $X^s$}\label{sect:Xs}

It is well known \citeaffixed[§4.8.1]{MR760629}{see, e.g.,} that the Chebyshev coefficients $a_m$ of $f(x)$ are given by 
the Fourier coefficients of $f(\cos\theta)$: 
\[
a_m = \frac{2}{\pi}\int_{-1}^1 \frac{f(x) T_m(x)}{\sqrt{1-x^2}}\,dx = \frac{2}{\pi}\int_0^\pi f(\cos\theta)\cos m\theta \,d\theta.
\]
Asymptotic analysis of Fourier integrals can now be used to determine the decay rate of the $a_m$: e.g., the function $f_s(x) = |x-\xi|^s$ with $-1<\xi<1$ and $s>0$ is of class $X^s$ since by the method of stationary phase \cite[§§3.11--3.13]{MR0435697}
\[
a_{m} = -\frac{4}{\pi}T_m(\xi)(1-\xi)^{s/2}\Gamma(1+s)\sin(\pi s/2) m^{-s-1} + o(m^{-s-1}) \qquad (m\to\infty).
\]
Alternatively (but often less sharp), decay estimates of Fourier coefficients based on
the smoothness properties of $f$ can be used; e.g.,
\cite[Thms. II.4.12]{MR0236587}:

\begin{quote}\emph{Let $f$ be defined on $[-1,1]$. If $f(\cos t)$ is $k-1$-times differentiable 
with a piecewise $k$-th derivative of bounded variation, then $f \in X^{k}$.}
\end{quote}
Since all derivatives of $\cos t$ exist and are bounded by the constant $1$, the smoothness properties
of $f(\cos t)$ can conveniently be inferred from those of $f(x)$ (but not vice versa). In particular, if $f$ itself is
 $k-1$-times differentiable with a piecewise $k$-th derivative of bounded variation, we still get $f \in X^{k}$. 

\subsubsection*{Remark} Denoting the total variation of that piecewise $k$-th derivative of $f$ by $V$, \citeasnoun[Thm.~7.1]{Trefethen} proved the explicit bound\footnote{We use Knuth's
notation of the $n$-th falling factorial power: $a^{\underline{n}}=a(a-1)\cdots (a-n+1)$.}
\[
|a_m| \leq \frac{2V}{\pi m^{\underline{k+1}}}\qquad (m\geq k+1);
\]
using it, \citeasnoun{Xiang} rendered the rate estimate (\ref{eq:claim}) in the explicit form
\[
|E_n(f)| \leq \frac{\pi V}{2 n^{\underline{k+1}}}
\]
if $n$ is sufficiently large (and, for Gauss quadrature, $k\geq 2$); an estimate that would asymptotically be, 
for $f(x)=|x|$, just
a factor of $2$ off the true state of affairs.

\section{Convergence rate of Clenshaw--Curtis quadrature}\label{sect:cc}

Clenshaw--Curtis quadrature on $[-1,1]$ is the interpolatory $n$-point quadrature rule that is derived from the nodes
\begin{equation}\label{eq:ccnodes}
x_k = \cos\left(\frac{k-1}{n-1}\pi\right)\qquad (k=1,\ldots,n).
\end{equation}
Now, it is well known that from $T_m(x) = \cos(m \arccos x)$ one reads off \emph{aliasing} due to undersampling, that is, with\footnote{Note that we do not need, for both quadrature rules studied in this paper, to consider \emph{odd} numbered Chebyshev polynomials: all their integrals and quadrature errors \emph{vanish} because of symmetry.} $m=2j(n-1)+2r$ and $-(n-2)\leq 2r \leq n-1$
\[
T_m(x_k) = T_{2|r|}(x_k);
\]
which implies, since Clenshaw--Curtis is exact for polynomials of degree $n-1$,
\[
I^C_n(T_m) = I^C_n(T_{2|r|}) = I(T_{2|r|}).
\]
Here, $I_n^C(f)$ denotes the quadrature formula as applied to $f$ and $I(f)$ the integral.
Therefore, as $m\geq n \to \infty$, the quadrature error $E_n^C(T_m)$ satisfies
\[
E_n^C(T_m) = I(T_m) - I(T_{2|r|}) = - \dfrac{2}{m^2-1} + \dfrac{2}{4r^2-1} = \dfrac{2}{4r^2-1} + O(n^{-2}).
\]
With $f\in X^s$, that is, $a_m=O(m^{-s-1})$ for some $s > 0$, we follow the ideas of \citeasnoun[p.~347]{MR0305555}
in estimating
\[
|E_n^C(f)| \leq \sum_{q=n}^\infty |a_{q}|\cdot |E_n^C(T_{q})| = O(S_1) + O(S_2),
\]
where 
\[
S_1 = \sum_{j=1}^\infty \sum_{|2r|< n} \dfrac{1/|4r^2-1|}{(2j(n-1)+2r)^{s+1}},\qquad
S_2 = n^{-2}\sum_{q=n}^\infty \frac{1}{q^{s+1}} = O(n^{-s-2}).
\]
Because of
\begin{equation}\label{eq:bounded}
\sum_{r=-\infty}^\infty \dfrac{1}{|4r^2-1|} = 2,\qquad \sum_{j=1}^\infty \frac{1}{j^{s+1}}= \zeta(s+1),
\end{equation}
we immediately see that $S_1 = O(n^{-s-1})$; hence
we obtain the rate estimate 
\begin{equation}\label{eq:rateC}
E_n^C(f) = O(n^{-s-1})\qquad (s>0),
\end{equation}
which proves the theorem of §\ref{sect:intro} in the Clenshaw--Curtis case.

\section{Convergence rate of Gauss quadrature I}\label{sect:GI}

As substitute for (\ref{eq:ccnodes}) there are \emph{asymptotic} formulas for the nodes $x_k$ of $n$-point Gauss quadrature (the zeros of the Legendre polynomial of degree $n$): a classical one of \citeasnoun{MR0093610} 
is, writing $\phi_k=(4k-1)\pi/(4n+2)$ for short,\footnote{\citeasnoun[p.~208]{MR0298934} stated this result 
with $O(n^{-3})$ instead of $O(k^{-2}n^{-1})$---citing as source \citeasnoun[p.~787]{MR0167642}, who had however \emph{misstated} 
the result of \citeasnoun{MR0093610}: Gatteschi's term $O(k^{-2}n^{-1})$ reduces to $O(n^{-3})$ 
only for those nodes $x_k$ that belong to a fixed interval in the interior of $[-1,1]$. However,
the calculations of \citeasnoun{MR0298934} are fairly easy to fix: in the end, their estimate of $E_n^G(T_m)$ turns out to be not affected at all.} 
\begin{equation}\label{eq:gaussnodes}
x_k = \cos\left(\phi_k + \tfrac18\cot(\phi_k)n^{-2} + O(k^{-2}n^{-1})\right)\qquad (1\leq k \leq n/2).
\end{equation}
Using this and an $O(n^{-1})$ bound on the weights, \citeasnoun[p.~211]{MR0298934} 
proved that the error in integrating the Chebyshev polynomials is\footnote{\citeasnoun{MR0298934} stated the remainder in the form $O(1/n)+O(\log n/n)$ for $m=O(n)$;
the explicit dependence on $m$ given here follows from noting that the quantities $h_i$ of their paper scale with $m/n$: the first remainder term estimates a weighted sum of $h_i^2$, the second a weighted sum of $|h_i|$.}
\[
E_n^G(T_m) = 
\begin{cases}
(-1)^j \dfrac{2}{4r^2-1} + O(m^2/n^3) + O(m\log n/n^2) & \quad -n< r <n,\\*[6mm]
(-1)^j \dfrac{\pi}{2} + O(m^2/n^3) + O(m\log n/n^2) & \quad r=\pm n,
\end{cases}
\]
if $2n\leq m=j(4n+2)+2r$ with $-n \leq r \leq n$ and $j\geq 0$. 
This way, aliasing holds asymptotically for $m=o(n^{3/2})$ only; for larger $m$, phase errors of
order~$O(1)$ will render the estimate useless. Still, because of $|T_m|\leq 1$ on $[-1,1]$ 
we get the uniform bound $|E_n^G(T_m)| \leq 4$. 
We now estimate $E_n^G(f) = E_n' + E_n''$
by splitting the Chebyshev expansion at an index of the order $O(n^{1+\epsilon})$ with some
$0<\epsilon < 1$ to be chosen later. Using the uniform
bound of $E_n^G(T_m)$ we thus get the tail estimate
\[
E_{n}''=\sum_{q=n^{1+\epsilon}}^\infty |a_{2q}|\cdot |E_n^G(T_{2q})| = O\left(\sum_{q=n^{1+\epsilon}}^\infty \dfrac{1}{q^{s+1}}\right) = O(n^{1-s\epsilon} n^{-s-1}).
\]
We are left with estimating the first $O(n^{1+\epsilon})$ terms of the Chebyshev expansion:
\[
E_n' = \sum_{q=n}^{n^{1+\epsilon}} |a_{2q}|\cdot |E_n^G(T_{2q})| = O(S_1') + O(S_2'),
\]
where
\begin{align*}
S_1' &= \sum_{j=1}^\infty \sum_{|r|<n} \dfrac{1/|4r^2-1|}{(j(4n+2)+2r)^{s+1}}
+ \sum_{j=1}^\infty \sum_{r=\pm n} \dfrac{1}{(j(4n+2)+2r)^{s+1}} + \dfrac{1}{n^{s+1}},\\*[2mm]
S_2' &= \dfrac{1}{n^3}\sum_{q=n}^{n^{1+\epsilon}} q^{1-s} + \frac{\log n}{n^2}\sum_{q=n}^{n^{1+\epsilon}} q^{-s}.
\end{align*}
From (\ref{eq:bounded}) we immediately see that $S_1' = O(n^{-s-1})$.
Likewise, we obtain 
\[
n^{s+1} S_2' =
\begin{cases}
O(n^{(2-s)\epsilon})&\quad 0<s<2,\\*[2mm]
O(\log n)&\quad s\geq 2.
\end{cases}
\]
Summarizing, the optimized choice $\epsilon = 1/2$ results in the rate estimate
\begin{equation}\label{eq:EnGfirst}
E_n^G(f) =
\begin{cases}
O(n^{-3s/2})      &\quad 0<s<2,\\*[2mm]
O(n^{-s-1}\log n)&\quad s\geq 2.
\end{cases}
\end{equation}
which proves the theorem of §\ref{sect:intro} in the Gauss case up to a factor $\log n$.

\section{Convergence rate of Gauss quadrature II}\label{sect:GII}

\citeasnoun{Xiang} observed that we can get rid of the logarithmic factor in (\ref{eq:EnGfirst}) by using
a refined estimate of
\citeasnoun[Thm.~1 and p.~199]{MR1345417}: upon replacing the
bound in (\ref{eq:gaussnodes}) by a later, sharper one also due to \citeasnoun{MR911648},\footnote{Luigi Gatteschi (1923--2007) worked for nearly 60 
years on the asymptotics of the zeros of special functions with a focus on explicit, useful error bounds; see \citeasnoun{MR2457088}.} namely
\begin{equation}\label{eq:gaussnodes}
x_k = \cos\left(\phi_k + \tfrac12\cot(\phi_k)(2n+1)^{-2} + O(k^{-3}n^{-1})\right)\qquad (1\leq k \leq n/2),
\end{equation}
and by using some improved, individual estimates of the weights, Petras proved, within the range $m = O(n^2)$, that
\[
|E_n^G(T_m)| = 
\begin{cases}
\dfrac{2+ O(m r/n^{2})}{|4r^2-1|} + O(m^4/n^{6}) + O(m^2 \log(n)/n^{4})& \, |r| < n,\\*[5mm]
\dfrac{\pi}{2} + O(m/n^{2}) + O(m^4/n^{6}) + O(m^2 \log(n)/n^{4}) & \, |r|=n,
\end{cases}
\]
where $2n\leq m= j(4n+2)+2r$ with $|r| \leq n$ and $0 \leq j = O(n)$. Thus, we obtain
\[
E_n' = \sum_{q=n}^{n^{1+\epsilon}} |a_{2q}|\cdot |E_n^G(T_{2q})| =O(S_1') + O(\tilde{S}_1') + O(\tilde{S}_2'),
\]
where $S_1' = O(n^{-s-1})$ is defined as in §\ref{sect:GI} and 
\begin{align*}
\tilde{S}_1' &= \sum_{j=1}^{n^\epsilon} \sum_{|r|<n} \dfrac{j r/|4r^2-1|/n}{(j(4n+2)+2r)^{s+1}}
+ \sum_{j=1}^{n^\epsilon} \sum_{r=\pm n} \dfrac{j/n}{(j(4n+2)+2r)^{s+1}} + \dfrac{1}{n^{s+2}},\\*[2mm]
\tilde{S}_2' &= \dfrac{1}{n^6}\sum_{q=n}^{n^{1+\epsilon}} q^{3-s} + \dfrac{\log n}{n^4}\sum_{q=n}^{n^{1+\epsilon}} q^{1-s}.
\end{align*}
By
\[
\sum_{r=-n}^n \dfrac{r}{|4r^2-1|} = O(\log n), \qquad \dfrac{1}{n}\sum_{j=1}^{n^\epsilon} j^{-s} = O(n^{\epsilon-1}),
\]
and, for $0<\epsilon<1$,  $O(n^{\epsilon-1}\log n ) = o(1)$ we get $\tilde{S}_1' = O(n^{-s-1})$.
Likewise
\[
n^{s+1} \tilde{S}_2' =
\begin{cases}
O(n^{(4-s)\epsilon}/n)&\quad 0 < s < 4,\\*[1mm]
O(\log n/n)&\quad s\geq 4.
\end{cases}
\]
Summarizing, though the optimal choice $\epsilon = 1/2$ just reproduces (\ref{eq:EnGfirst})
for $0<s<2$, it results, this time, in the rate estimate
\begin{equation}
E_n^G(f) = O(n^{-s-1}) \qquad (s\geq 2),
\end{equation}
which finally proves the Gauss case of the theorem of §\ref{sect:intro}.

\section{Open problems}
We leave the following open problems as challenges to the reader; their solution would require further, significant technical refinements of the methods used in this paper: to prove that, for $f\in X^s$,
\begin{itemize}
\item the convergence rate is $O(n^{-s-1})$ for Gauss quadrature if $0<s<2$;\\*[-3.5mm]
\item $|E_{n}^G(f)/E_n^C(f)|$ and its reciprocal stay uniformly bounded (cf. Fig.~\ref{fig:1}).
\end{itemize}

\subsubsection*{Acknowledgements} The authors thank Nick Trefethen for his continuing interest in this work and for his comments on some preliminary versions of the manuscript.

\bibliographystyle{kluwer}
\bibliography{quadrature}
\nopagebreak
\end{document}